\documentclass{article}

\usepackage{amsmath, amssymb, amsthm} 

\usepackage{hyperref} 

\renewcommand{\a}{\alpha} 
\renewcommand{\b}{\beta} 

\newcommand{\Gr}{\mathrm{Gr}} 
\newcommand{\St}{\mathrm{St}} 
\newcommand{\pr}{\mathrm{pr}} 
\newcommand{\Ad}{\mathrm{Ad}} 
\newcommand{\DET}{\mathrm{Det}} 
 
\newcommand{\res}{\mathrm{res}}
\newcommand{\Det}{\mathrm{Det}} 
\newcommand{\tr}{\mathrm{tr}} 
\renewcommand{\top}{\mathrm{top}}
\renewcommand{\det}{\mathrm{det}} 
\newcommand{\Map}{\mathrm{Map}} 
\newcommand{\Lie}{\mathrm{Lie}}

\newcommand{\cB}{\mathcal{B}} 
 
\newcommand{\cE}{\mathcal{E}} 
\newcommand{\cR}{\mathcal{R}}

\newcommand{\CC}{\mathbb{C}}

\newcommand{\g}{\mathfrak{g}}

\theoremstyle{plain}
\newtheorem{theorem}{Theorem}
\newtheorem{lemma}[theorem]{Lemma}
\newtheorem{proposition}[theorem]{Proposition}

\theoremstyle{definition}

\theoremstyle{remark}

\title{Geometry of infinite dimensional Grassmannians and the Mickelsson-Rajeev cocycle}
\author{
Danny Stevenson \thanks{Fachbereich Mathematik,
Universit{\"a}t Hamburg, Hamburg, 20146, Germany,  \newline
 Email: \texttt{stevenson@math.uni-hamburg.de}} 
}
\begin{document} 
\maketitle 
\begin{abstract}
In their study of the representation theory of loop groups, Pressley and Segal introduced a determinant line bundle 
over an infinite dimensional Grassmann manifold.  Mickelsson and Rajeev subsequently generalized 
the work of Pressely and Segal to obtain representations of the groups $\Map(M,G)$ where $M$ is an 
odd dimensional spin manifold.  In the course of their work, Mickelsson and Rajeev introduced for any $p\geq 1$, 
an infinite dimensional Grassmannian $\Gr_p$ and a determinant line bundle $\Det_p$ over it, generalizing the constructions 
of Pressley and Segal.  The definition of the line bundle $\Det_p$ requires the notion of a regularized determinant for bounded operators.  
In this note we specialize to the case when $p=2$ (which is relevant for the case when $\dim M = 3$) and consider 
the geometry of the determinant line bundle $\Det_2$.  We construct explicitly a connection on $\Det_2$ and 
give a simple formula for its curvature.  From our results we obtain a geometric derivation of the Mickelsson-Rajeev 
cocycle.   
\end{abstract}

\section{Introduction} 

In the paper \cite{MR}, the authors construct representations of the groups 
$\Map(M,G)$, generalizing the methods of Pressley 
and Segal \cite{PS} for constructing representations of loop groups.  Here $M$ is a 
compact spin manifold of odd dimension and $G$ is a compact Lie group.  
In the work of Pressley and Segal a fundamental role was played by the \emph{restricted 
general linear group} $GL_{\res}$ and the \emph{restricted Grassmannian} $\Gr_{\res}$ 
associated to a polarized Hilbert space $H = H_+\oplus H_-$.  $GL_{\res}$ and $\Gr_{\res}$ were defined 
relative to a certain Schatten ideal, namely the Hilbert-Schmidt operators.  Recall 
that for any $p\geq 1$ one can define ideals $L_p$ --- the {\em Schatten ideals} --- in the space $\cB(H)$ of bounded 
operators on $H$ (see for example \cite{Simon}).  When the Schatten index $p=1$, the ideal $L_1$ is 
just the ideal of trace class operators on $H$, and when $p=2$ the ideal $L_2$ is 
the ideal of Hilbert-Schmidt operators on $H$, as we have mentioned.  These 
Schatten ideals play an important role in non-commutative geometry \cite{Connes}.  
They arise also in the work of Mickelsson and Rajeev.  In \cite{MR} 
the group $GL_{\res}$ was generalized to the group $GL_{(p)}$, for any Schatten index $p\geq 1$.  
The group $GL_{(p)}$ (not to be confused with the general linear group of a $p$-dimensional 
vector space!) was defined to be the group of invertible 
operators $g$ on $H$ such that with respect to the polarization $H = H_+\oplus H_-$ 
the operator $g$ has $2\times 2$ block operator form 
$$ 
\left(\begin{matrix} a & b \\ c & d \end{matrix} \right) 
$$ 
in which the off diagonal blocks $b$ and $c$ belong to the Schatten class $L_{2p}$.  
$\Gr_p$ is the associated Grassmannian, again defined relative to $L_{2p}$.  We shall recall 
the definition of $\Gr_p$ in greater detail in Section~\ref{sec2}.  In the 
framework of Mickelsson and Rajeev $GL_{\res}$ and $\Gr_{\res}$ correspond to $GL_{(1)}$ and $\Gr_1$ respectively.  
The Schatten index $p$ arises in \cite{MR} in the following way.  
In that paper it is shown that there is an embedding of the group $\Map(M,G)$ into the general linear 
group $GL_{(p)}$ of a certain polarized Hilbert space provided that $p$ exceeds a certain 
bound, related to the dimension of $M$.  This generalizes the embedding defined by Pressley and Segal of the loop group 
$\Map(S^1,G)$ into $GL_{\res}$.

In the case of the ordinary Grassmannian $\Gr(V)$ associated to a finite dimensional 
vector space $V$, there is a canonical holomorphic determinant line bundle $\Det$ defined over 
$\Gr(V)$.  If $W$ is a subspace of $V$ belonging to some connected component of 
$\Gr(V)$ then the fiber of $\Det$ at $W$ is the top exterior power $\Lambda^\top W$.  
In the case of the infinite dimensional Grassmannian $\Gr_1$ the notion of the 
top exterior power loses its meaning.  Nevertheless, Pressley and Segal construct a well defined holomorphic determinant line 
bundle $\Det$ on $\Gr_{\res}$, and moreover show that there is a central extension of groups 
\begin{equation} 
\label{eq: PS extn} 
1\to \CC^*\to \widehat{GL}_{\res}\to GL_{\res}\to 1 
\end{equation}
with the property that the group $\widehat{GL}_{\res}$ acts on the space $\Gamma$ of holomorphic sections of the dual 
bundle $\Det^*$.  Pressley and Segal show that $\Gamma$ can be interpreted as the fermionic 
Fock space construction on $H$.  Inside $GL_{\res}$ is the subgroup $U_{\res}$ consisting of 
all unitary operators in $GL_{\res}$.  Corresponding to $U_{\res}$ is a subgroup $\widehat{U}_{\res}$ 
of $\widehat{GL}_{\res}$ and it turns out that the irreducible representation of $\widehat{GL}_{\res}$ on 
$\Gamma$ restricts to an irreducible unitary 
representation of $\widehat{U}_{\res}$.  This irreducible representation of $\widehat{U}_{\res}$ is used to construct the 
`basic' positive energy representation of the loop group $LU_n$ using an embedding $LU_n\subset U_{\res}$.  
This is the construction the authors in \cite{MR} worked towards generalizing to the groups $\Map(M,G)$.  

Using a notion of regularized determinant \cite{Simon} for invertible operators in $1 + L_p$ Mickelsson and 
Rajeev construct \emph{regularized determinant line bundles} $\Det_p$ on $\Gr_p$.  The central 
extension~\eqref{eq: PS extn} above is replaced by an extension 
\begin{equation} 
\label{eq: MR extn} 
1\to \Map(\Gr_p,\CC^*)\to \widehat{GL}_{(p)}\to GL_{(p)}\to 1 
\end{equation}
which is now non-central.  An extension of $\widehat{GL}_{(p)}$ is obtained on the space of smooth 
sections of $\Det_p^*$.  This extension satisfied a positive energy condition, however it was later shown \cite{Mic2, Pickrell} that this was not a unitary 
extension when $p=2$.  As mentioned above the appearance of the Schatten indices $p$ 
in the work of Mickelsson and Rajeev arises from an embedding $\Map(M,G)\subset 
GL_{(p)}$, where $p$ depends on the dimension of $M$.  The case when the dimension of $M$ is 3 
corresponds to $p=2$.    

The purpose of this note is to study the geometry of the Grassmannian $\Gr_2$.  
Our main result (Proposition~\ref{prop: connection}) gives an explicit construction of a 
connection on $\Det_2$ and an explicit and simple formula for the corresponding 
curvature 2-form.  We point out that this is not as trivial as it seems.  It is easy to 
do this for the determinant line bundle over $\Gr_{\res}$ studied by Pressley and 
Segal, but the case of $\Det_2$ is much more delicate.  The same difficulties arise in 
finding closed formulas for certain universal Schwinger cocycles: when $p=1$  
there is the well known Kac-Peterson cocycle (see \cite{PS} Proposition 6.6.5), when $p=2$ 
 a considerably more difficult calculation produces the 
Mickelsson-Rajeev cocycle (see \cite{MR} and the discussion below).  No closed formula 
is known for these universal cocycles for arbitrary $p$ (but see \cite{FT} for a conjectural formula).     

We should also mention that the line bundles $\Det_p$ are not just holomorphic line bundles, 
but they also carry a Hermitian structure as well.  Therefore there is a canonical connection 
on each $\Det_p$ compatible with both the Hermitian and holomorphic structures.  There is a 
formula for the curvature of this canonical connection however it appears to be quite difficult 
to derive a simple expression for it, at least for $p>1$.  This canonical connection 
featured also in Quillen's paper \cite{Q2} however there he managed to identify its curvature with 
a certain K\"{a}hler form.

One consequence of our result is that we obtain a simple and explicit formula for 
a de Rham representative of the first Chern class $c_1(\Gr_2)$ of $\Gr_2$.  Recall 
that the topology of the Grassmannians $\Gr_p$ is well understood \cite{Palais, Q1}.  For 
any $p<q$ there is a natural inclusion $\Gr_p\subset \Gr_q$ and this turns out to be a 
homotopy equivalence.  The infinite dimensional manifolds   
$\Gr_p$ give smooth models for the classifying space 
of even $K$-theory $K^0$ and in fact it turns out that the de Rham theorem holds for them.  
This was exploited by Quillen in \cite{Q1} where he gave explicit formulas in terms of 
contour integrals for differential form representatives of the Chern classes $c_n(\Gr_p)$ for all $n$ and $p$. The expressions 
he obtained however were not easy to evaluate directly.  

Another by-product of our construction of a connection on $\Det_2$ is that we are able 
to give a geometric derivation of the \emph{Mickelsson-Rajeev cocycle} associated to 
the extension of Lie algebras 
\begin{equation} 
\label{eq: inf MR extn} 
0\to \Map(\Gr_2,\CC)\to \widehat{\mathfrak{gl}}_{(2)}\to \mathfrak{gl}_{(2)}\to 0 
\end{equation} 
which is the infinitesimal version of~\eqref{eq: MR extn} (again note the potential 
confusion with the finite dimensional Lie algebra $\mathfrak{gl}(2)$!).  
This cocycle was derived in \cite{MR} however the computation 
used there was rather involved.  An alternative, algebraic derivation 
can be found in \cite{Langmann}.   Using the curvature 2-form of our connection on $\Det_2$ 
we give an alternative expression for the cocycle associated to the extension~\eqref{eq: inf MR extn}.  
This does not agree on the nose with the cocycle obtained by Mickelsson and Rajeev, but we give an explicit formula for a 
coboundary relating the two cocycles.   

The extension~\eqref{eq: MR extn} is an example of an extension of groups that arises naturally 
whenever one has a line bundle $L$ on a manifold $M$ on which a Lie group $G$ acts.  
In this situation there is a canonical extension 
\begin{equation} 
\label{eq: can extn}
1\to \Map(M,\CC^*)\to \widehat{G}\to G\to 1 
\end{equation}  
of $G$ by the abelian group $\Map(M,\CC^*)$.  
The Mickelsson-Rajeev extension is a special case of this canonical extension.  We describe 
two methods for associating a Lie algebra 2-cocycle to the extension of Lie algebras associated to~\eqref{eq: can extn}.  
The first method is geometric, using a connection on $L$.  The second method is more algebraic, 
requiring a knowledge of the local structure of~\eqref{eq: can extn}.  The two methods lead to different 
2-cocycles in general.  We give a formula, which appears to be new, for a coboundary 
relating these 2-cocycles.  
  
In summary then this paper is as follows. In Section~\ref{sec1} we discuss the two methods for 
associating a cocycle to the extension of Lie algebras associated to~\eqref{eq: can extn} and 
derive a formula for a coboundary relating the two different cocycles we obtain.    In Sections~\ref{sec2} and~\ref{sec3} we provide some background on the infinite 
dimensional Grassmann manifold $\Gr_2$ and the determinant line bundle $\Det_2$.  
In Section~\ref{sec4} we construct a connection 1-form on $\Det_2$ and compute its 
curvature.  Section~\ref{sec5} contains a comparison of the geometric cocycle describing the extension~\eqref{eq: inf MR extn} 
that we obtain from the curvature with the Mickelsson-Rajeev cocycle.  Two slightly complicated calculations are 
contained in the appendices.     

\section{General remarks on Lie algebra $2$-cocycles} 
\label{sec1} 

Suppose that $M$ is a $\g$-module for a Lie algebra $\g$.  
We can consider $M$ as an abelian Lie algebra and 
consider extensions of Lie algebras
$$ 
0\to M\to \widehat{\g}\stackrel{p}{\to} \g\to 0. 
$$ 
If the linear map underlying the homomorphism 
$p$ admits a section (so that $\widehat{\g}\cong \g\oplus M$ as a vector space) then one can associate to 
the extension a Lie algebra $2$-cocycle $\omega$ with values in 
the $\g$-module $M$; thus $\omega$ is a linear map 
$\omega\colon \Lambda^2 \g\to M$ 
such that 
\begin{multline} 
\label{eq: cocycle condition} 
\omega([\xi,\eta],\zeta) - \omega([\xi,\zeta],\eta) + \omega([\eta,\zeta],\xi) \\ 
- \xi\cdot\omega(\eta,\zeta) + \eta\cdot\omega(\xi,\zeta) - \zeta\cdot\omega(\xi,\eta) = 0 
\end{multline} 
Conversely, given such a cocycle one may use it to 
twist the Lie bracket on $\g\oplus M$ to obtain a new Lie algebra $\widehat{\g}$ 
fitting into an extension of Lie algebras as above. This is a brief summary of the well-known 
theorem that isomorphism classes of extensions of $\g$ by the abelian 
Lie algebra $M$ are classified by the Lie algebra cohomology 
group $H^2(\g,M)$.  

\bigskip 

A nice example of such an extension of Lie algebras arises as the infinitesimal version 
of the extension of Lie groups~\eqref{eq: can extn} mentioned in the Introduction.  
Suppose $L$ is a line bundle on a manifold $M$ on which 
a Lie group $G$ acts.  Then there is canonically associated 
to $G$ an extension of Lie groups 
\begin{equation} 
\label{eq: group extension} 
1\to \Map(M,\CC^*)\to \widehat{G}\to G\to 1 
\end{equation}
where $\widehat{G}$ is the subgroup of the group of 
bundle automorphisms of $L$ consisting of automorphisms which cover the action of 
$G$ on $M$.    
There is a corresponding infinitesimal version 
of this; if $\mathfrak{g}$ denotes the Lie algebra of 
$G$ then we have the extension of Lie algebras 
\begin{equation} 
\label{eq: lie alg extension}
0\to \Map(M,\CC)\to \widehat{\mathfrak{g}}\to \mathfrak{g}\to 0.   
\end{equation}  
Here a vector $\xi\in \mathfrak{g}$ acts on a function $f\in \Map(M,\CC)$ by 
\begin{equation} 
\label{eq: inf action} 
(\xi\cdot f)(x) = \frac{d}{dt}\Big|_{t=0}f(x\exp(t\xi)). 
\end{equation}
In fact every vector $\xi\in \mathfrak{g}$ generates a vector field $\hat{\xi}$ on $M$, the fundamental 
vector field generated by the infinitesimal action of $\xi$.  The value of $\hat{\xi}$ on a function 
$f$ on $M$ is given by precisely the same derivative formula as in~\eqref{eq: inf action}.  

It is important to note that in general the extensions~\eqref{eq: group extension} 
and~\eqref{eq: lie alg extension} are not central extensions.  We will call such 
an extension of groups by an abelian normal subgroup an \emph{abelian extension} of groups, and we will also call an extension of 
Lie algebras by an abelian ideal, an \emph{abelian extension} of Lie algebras.   

By the classification theorem mentioned above, 
the abelian extension of Lie algebras~\eqref{eq: lie alg extension} will be described by a Lie 
algebra $2$-cocycle $\omega$ on $\mathfrak{g}$ 
with values in the $\mathfrak{g}$-module $\Map(M,\CC)$.  
Since the extension of groups, and hence the associated extension~\eqref{eq: lie alg extension},  
is completely determined by the line bundle $L$ on $M$ 
and the action of the group $G$ on $M$, then one should expect  
to find a formula for the cocycle $\omega$ in terms of some 
geometric data on $L$.  
Indeed, 
if $L$ comes equipped with a connection $\nabla$ whose 
curvature $2$-form is $F_{\nabla}$ then one can describe 
$\omega$ as follows: 
\begin{equation} 
\label{eq: lie alg 2 cocycle} 
\omega(\xi,\eta) = -F_{\nabla}(\hat{\xi},\hat{\eta})
\end{equation}
where $\hat{\xi}$ and $\hat{\eta}$ are the fundamental vector fields on 
$M$ generated by the infinitesimal action of $\xi,\eta\in \mathfrak{g}$.  
The condition~\eqref{eq: cocycle condition} that $\omega$ 
is a Lie algebra $2$-cocycle is exactly the condition that the 
curvature $F_{\nabla}$ is a closed $2$-form on $M$.  

Another method to compute the cocycle $\omega$ is to use  
the local structure of the group $\widehat{G}$.  Here it is 
important to realise that \emph{locally} the underlying 
manifold of the Lie group $\widehat{G}$ is 
a product of $G$ and $\Map(M,\CC^*)$ but this is 
not in general true \emph{globally}.  In other words 
$\widehat{G}$ is a locally trivial principal $\Map(M,\CC^*)$ 
bundle over $G$.  
In this method one chooses a local section $\sigma$ defined in a neighbourhood of $1$ of the 
map underlying the homomorphism $\widehat{G}\to G$  
and defines a 2-cocycle $\omega(\xi,\eta)$ by   
\begin{equation}
\label{eq: alg 2-cocycle} 
\frac{\partial^2}{\partial s\partial t}\Big|_{s=t=0} 
\sigma(e^{t\xi})\sigma(e^{s\eta})\sigma(e^{-t\xi})\sigma(e^{-s\eta}) = ([\xi,\eta],\omega(\xi,\eta))  
\end{equation}
where $\exp(t\xi)$ and 
$\exp(s\eta)$ are $1$-parameter subgroups.   
Here the right hand side should be understood in terms of the 
splitting $\widehat{\mathfrak{g}}\cong \mathfrak{g}\oplus 
\Map(M,\CC)$ of the exact 
sequence~\eqref{eq: lie alg extension} of Lie algebras defined by $d\sigma$, 
the derivative of the local section $\sigma$ at the identity.  
Thus another way to think of $\omega(\xi,\eta)$ is as the 
familiar expression 
$$ 
\omega(\xi,\eta) = [d\sigma(\xi),d\sigma(\eta)] - d\sigma[\xi,\eta].
$$   
We remark that the expression~\eqref{eq: alg 2-cocycle} can be difficult to evaluate, cf. equation (4.12) 
in \cite{MR}.  

We have then two different descriptions of the 
Lie algebra $2$-cocycle associated to the extension~\eqref{eq: lie alg extension} --- let us denote the 
$2$-cocycle obtained from the geometric data (i.e the curvature) 
by $\omega_G(\xi,\eta)$ and the $2$-cocycle 
obtained from the local, algebraic structure by $\omega_A(\xi,\eta)$.  
Since the cocycles $\omega_G(\xi,\eta)$ and $\omega_A(\xi,\eta)$ 
define the same extension of Lie algebras they should be cohomologous. In fact 
one can write down an explicit formula for a coboundary\footnote{We haven't been able to find a reference in which the 
formula above for the coboundary is described, but we would surprised 
if it were not known. } relating 
them, this is the content of our first proposition.  

Before we state the proposition, we will make a remark about line bundles.  
Throughout the paper we will blur the distinction between line bundles and 
principal $\CC^*$ bundles.  To every line bundle $L$ is associated a principal 
$\CC^*$ bundle $L^+$, its principal frame bundle.  The association of $L^+$ to 
$L$ sets up an equivalence of categories between the category of line bundles on $M$ 
and the category of principal $\CC^*$ bundles on $M$.  It is well known (see for example 
\cite{Brylinski}) that this equivalence between line bundles and principal $\CC^*$ bundles 
extends to connections: to every connection $\nabla$ on a line bundle $L$ there is 
associated a connection 1-form $A$ on the principal $\CC^*$ bundle $L^+$ and conversely.  
For more details we refer to \cite{Brylinski}.  With these remarks made, we will freely pass 
between line bundles and principal $\CC^*$ bundles without further comment.  

We recall 
that $L$ was a line bundle over the $G$-manifold $M$ equipped with a connection $\nabla$.  
We will denote by $A$ the corresponding connection 1-form and we will denote by $\sigma$ a 
local section of $\widehat{G}\to G$ defined in a neighborhood of the identity.   
\begin{proposition} 
\label{prop: coboundary}
Let $L$, $M$, $G$, $A$ and $\sigma$ be as above.  Then the two Lie 
algebra cocycles $\omega_G(\xi,\eta)$ and $\omega_A(\xi,\eta)$ are 
related by the coboundary $b(\xi)$ defined by 
\begin{equation}
\label{eq: coboundary formula} 
b(\xi) = A(\widehat{d\sigma(\xi)}).  
\end{equation} 
In other words we 
have 
$$ 
\omega_A(\xi,\eta) = \omega_G(\xi,\eta) + \xi\cdot b(\eta) - \eta\cdot b(\xi) - b([\xi,\eta]). 
$$ 
\end{proposition} 
 First 
suppose that $\xi$ is a vector in $\mathfrak{g}$.  Then 
$d\sigma(\xi)$ is a vector in $\widehat{\mathfrak{g}}$ and 
we can consider the fundamental vector field $\widehat{d\sigma(\xi)}$ 
on $L$ induced by the infinitesimal action of $d\sigma(\xi)$.  We 
can also consider the horizontal lift $(\hat{\xi})_H$ of 
the fundamental vector field $\hat{\xi}$ on $M$.       
The vector field $\widehat{d\sigma(X)} - (\hat{\xi})_H$ on 
$L$ is vertical with respect to the $\CC^*$ action on $L$.  
Consider 
\begin{multline*} 
[\widehat{d\sigma(\xi)},\widehat{d\sigma(\eta)}]  = 
[\widehat{d\sigma(\xi)} - (\hat{\xi})_H,\widehat{d\sigma(\eta)} - 
(\hat{\eta})_H] + [(\hat{\xi})_H,\widehat{d\sigma(\eta)} - (\hat{\eta})_H] \\
+ [\widehat{d\sigma(\xi)} - (\hat{\xi})_H,(\hat{\eta})_H] + [(\hat{\xi})_H,(\hat{\eta})_H] 
\end{multline*}
The term $[\widehat{d\sigma(\xi)} - (\hat{\xi})_H,\widehat{d\sigma(\eta)} - 
(\hat{\eta})_H]$ in this expression vanishes, since it is 
a bracket of the form $[\hat{\a},\hat{\b}]$ where $\a$ and $\b$ are complex numbers.  
We have, since $F_A(\hat{\xi},\hat{\eta}) = - A([(\hat{\xi})_H,(\hat{\eta})_H])$,  
$$ 
A([\widehat{d\sigma(\xi)},\widehat{d\sigma(\eta)}]) = 
A([\widehat{d\sigma(\xi)} - (\hat{\xi})_H,(\hat{\eta})_H]) + 
A([(\hat{\xi})_H,\widehat{d\sigma(\eta)} - (\hat{\eta})_H]) - F_A(\hat{\xi},\hat{\eta}). 
$$ 
If $V$ and $W$ are vectors in $\Lie(\widehat{G})$ then we 
have the relation 
$\widehat{[V,W]} = [\widehat{V},\widehat{W}]$ between 
fundamental vector fields.  It follows therefore that 
we have 
\begin{align*}
c(\xi,\eta)  & = A([\widehat{d\sigma(\xi)},\widehat{d\sigma(\eta)}]) - A(\widehat{d\sigma[\xi,\eta]}) \\ 
& =  -F_A(\hat{\xi},\hat{\eta}) + A([(\hat{\xi})_H,\widehat{d\sigma(\eta)} - (\hat{\eta})_H]) 
- A([(\hat{\eta})_H,\widehat{d\sigma(\xi)} - (\hat{\xi})_H]) \\
& \phantom{=  -F_A(\hat{\xi},\hat{\eta}) + A([(\hat{\xi})_H,\widehat{d\sigma(\eta)} - (\hat{\eta})_H]) 
- A([(\hat{\eta})_H, } - A(\widehat{d\sigma[\xi,\eta]})
\end{align*} 
To simplify the terms $A([(\hat{\xi})_H,\widehat{d\sigma(\eta)} - (\hat{\eta})_H])$ 
and $A([(\hat{\eta})_H,\widehat{d\sigma(\xi)} - (\hat{\xi})_H])$ consider 
more generally $A([X,Y])$ where the vector field $X$ is vertical and the 
vector field $Y$ is horizontal.  Then it is easy to see that we have $A([X,Y]) = -\mathcal{L}_YA(X)$.
Therefore we can write 
$$ 
A([(\hat{\xi})_H,\widehat{d\sigma(\eta)} - (\hat{\eta})_H]) = 
\mathcal{L}_{(\hat{\xi})_H}A(\widehat{d\sigma(\eta)} - (\hat{\eta})_H) = 
\mathcal{L}_{(\hat{\xi})_H}A(\widehat{d\sigma(\eta)}). 
$$ 
The vector field $\widehat{d\sigma(\eta)}$ on $L$ is invariant under the 
action of $\CC^*$: $dR_z(\widehat{d\sigma(\eta)}_p = \widehat{d\sigma(\eta)}_{pz}$, 
and so $A(\widehat{d\sigma(\eta)})$ descends to a function on $M$, 
which we will continue to denote $A(\widehat{d\sigma(\eta)})$.  
If $V$ is a vector field on $M$ and $f$ is a function on $M$ then we calculate the Lie derivative 
$\mathcal{L}_V(f)$ at $x\in M$ by choosing a path $\gamma\colon 
(-\epsilon,\epsilon)\to M$ through $x$ with $\gamma'(0) =V$ and taking 
the derivative 
$$ 
\frac{d}{dt}\Big|_{t=0}f(\gamma(t)). 
$$ 
If $\gamma_H$ denotes a horizontal lift of $\gamma$ then this derivative 
is also equal to 
$$ 
\frac{d}{dt}\Big|_{t=0}\hat{f}(\gamma_H(t)). 
$$ 
where $\hat{f}(p) = f(\pi(p))$ ($\pi\colon L\to M$ denotes the projection).  
The point  of this discussion is that we can identify $\mathcal{L}_{(\hat{\xi})_H}A(\widehat{d\sigma(\eta)})$ 
with 
$$ 
\mathcal{L}_{\hat{\xi}}A(\widehat{d\sigma(\eta)}) = \xi\cdot A(\widehat{d\sigma(\eta)}) 
$$ 
where $\xi\cdot f$ denotes the action of $\xi\in \mathfrak{g}$ on 
a function $f$ in the $\mathfrak{g}$-module 
$\Map(M,\CC)$.  From here it is easy to see that $b(\xi) = A(\widehat{d\sigma(\xi)})$ is the required coboundary.

For the remainder of the paper we would like to study the example 
of this general situation mentioned in the introduction; namely 
the extension of groups described in \cite{MR} associated to a 
the determinant line bundle $\Det_2$ over the infinite dimensional 
Grassmannian manifold $\Gr_2$.  In the sections that follow we briefly 
review the construction of this determinant line bundle 
and discuss its geometry.  
 
\section{The geometry of the Grassmannians $\Gr_p$} 
\label{sec2}

The description of the determinant line bundle $\Det_p$ requires a basic knowledge 
of the Schatten ideals $L_p$ in the algebra of bounded operators 
$\cB(H)$ on a separable complex Hilbert space $H$.  
Briefly $L_p$ is defined to be the set of all operators 
$A\in \cB(H)$ such that $\tr(AA^*)^{p/2}< \infty$.  
If $A\in L_p$ then we write $||A||_p = (\tr(AA^*)^{p/2})^{1/p}$. 
It can be shown that $||\cdot  ||_p$ defines 
a norm on $L_p$ and that $L_p$ is an ideal in $\cB(H)$ for 
any $p\geq 1$.  As remarked earlier in the introduction,     
$L_1$ consists of the trace class operators on 
$H$ and $L_2$ consists of the Hilbert-Schmidt operators.  
The ideals $L_p$ share many of the properties 
of the measure spaces $L^p(X)$; for example if $A\in L_{q}$, 
 $B\in L_{r}$ and $p^{-1} =  q^{-1} + r^{-1}$ 
then $AB\in L_p$ and $||AB||_p \leq ||A||_{q}||B||_{r}$.  
For more details the reader should consult \cite{Simon}. 

To describe the line bundles $\Det_p$ we first need to describe the 
manifolds $\Gr_p$ over which they are defined.  $\Gr_p$ is an infinite 
dimensional Grassmannian manifold associated to a complex, 
infinite dimensional, separable Hilbert space $H$ which is equipped 
with a polarization $H = H_+\oplus H_-$.  To this polarization we 
can associate the operator 
$$
\epsilon = \left(\begin{array}{rr} 
1 & 0 \\ 0 & - 1 \end{array}\right)
$$ 
which 
is $1$ on the subspace $H_+$ and $-1$ on the subspace 
$H_-$.  $\epsilon$ is self adjoint and $\epsilon^2 = 1$.  In terms of the projections $\pr_+$ and $\pr_-$ 
onto the subspaces $H_+$ and $H_-$ respectively, $\epsilon$ can be written 
as $\epsilon = 2\pr_+ - 1$.  This operator $\epsilon$ plays a useful role in defining 
the manifold $\Gr_p$ as we now explain.  

In the case of the ordinary Grassmannian $\Gr(V)$ associated to a vector space $V$, 
there are several different ways to describe points in $\Gr(V)$.  These different ways 
are described in \cite{Q1}.  One can either think of a point of $\Gr(V)$ as subspace $W\subset V$, 
or equivalently we can replace the subspace $W$ with the orthogonal projection $P_W$ onto it.  
We can replace the orthogonal projection $P_W$ with the self adjoint involution $F$ of $V$ 
defined by $F = 2P_W - 1$.  Clearly we can move back and forth between projections and involutions 
this way.  Finally, the group $GL(V)$ acts transitively on $\Gr(V)$, and this leads to another description 
of $\Gr(V)$ as a homogenous space.  These four descriptions of points in Grassmannians persist 
to the infinite dimensional case of $\Gr_p$.  From \cite{MR,Q1} we have the following descriptions 
of points in $\Gr_p$:  
\begin{enumerate} 
\item a point of $\Gr_p$ can be thought of as a subspace $W\subset H$ such that the orthogonal 
projections $\pr_+\colon W\to H_+$ and $\pr_-\colon W\to H_-$ are Fredholm and $L_{2p}$ operators respectively, 

\item a point of $\Gr_p$ can be thought of as a self adjoint projection $P$ on $H$ such that $[P,\epsilon]\in L_{2p}$, 

\item a point of $\Gr_p$ can be thought of as a self adjoint bounded operator $F$ on $H$ such that $F - \epsilon \in L_{p}$ 
and $F^2 = 1$,

\item $\Gr_p$ can be thought of as the homogenous space $GL_{(p)}/B_{(p)}$ where $B_{(p)}$ is the subgroup 
of $GL_{(p)}$ consisting of invertible operators with block diagonal decomposition of the form 
$$ 
\left(\begin{array}{cc} 
a & b \\ 0 & d \end{array}\right) 
$$  
\end{enumerate}
Each of these definitions have their own advantages, for instance from 4 it is clear that $\Gr_{(p)}$ has a natural structure as 
a complex Banach manifold.  We find the description in 3 in terms of involutions the most convenient for our purposes.  
For the remainder of this paper we shall only be interested in the $L_4$ Grassmannian $\Gr_2$ and we will take a moment 
to amplify the description in 3 for this case.  Points in the $L_4$ Grassmannian are self adjoint bounded involutions $F$ on $H$ 
such that $F - \epsilon \in L_4$.  With respect to the polarization $H = H_+ \oplus H_-$ we can write 
$$ 
F = \left(\begin{array}{cc} F_{11} & F_{12} \\ F_{21} & F_{22} \end{array}\right) 
$$ 
Thus $F_{11}$ and $F_{22}$ are self adjoint and $F_{21} = F_{12}^*$.  Since $F - \epsilon \in L_4$, we must have 
$(F - \epsilon)^2 \in L_2$.  From here we see that $F_{11} - 1\in L_2$ and $F_{22} + 1\in L_2$.  
Since $F^2 = 1$ and $F_{11} \in 1 + L_2$ we see also that $F_{12}F_{12}^*\in L_2$ 
and hence $F_{12}\in L_4$.  

Associated to each such involution $F$ is a subspace $W$ satisfying the conditions of 1 above.  In particular 
the orthogonal projection $\pr_+\colon W\to H_+$ is a Fredholm operator, the index of which is called 
the \emph{virtual dimension} of $W$ (see \cite{PS}).  
In general there are many 
components of $\Gr_2$ and in fact these components are labelled by the 
virtual dimension.  In this paper we will just be concerned 
with the connected component $(\Gr_2)_0$ of $\Gr_2$ consisting of planes 
$W$ of virtual dimension zero.  For this reason we will indulge in a slight abuse of notation and write 
$\Gr_2$ when we really mean $(\Gr_2)_0$.  

Over $\Gr_2$ there is defined (see \cite{Mic1, MR, PS}) 
a `Steifel' bundle $\St_2\to \Gr_2$ of `admissible frames'.  
This is a principal bundle with structure group $GL^2$, 
the invertible operators $g$ on $H_+$ such that 
$g-1$ is Hilbert-Schmidt.  The space
$\St_2$ can be described by  
$$ 
\St_2 = \{w\colon H_+\to H|\ w\ \text{is injective},\ \pr_+w - 1\in L_2, 
\pr_- w\in L_4\}. 
$$ 
If $w\in \St_2$  then the image $w(H_+)$ is a subspace $W$ of $H$ satisfying the 
conditions of 1 above, and hence $W$ defines a point of $\Gr_2$.  The invertible operator 
$w\colon H_+\to H$ defines a basis for $W$; this is an admissible frame in the sense of 
\cite{Mic1, PS}.  
Thus $\St_2$ is an open subset of the Banach space 
which is the subspace of $\cB(H_+,H)$ consisting of all 
bounded operators $T\colon H_+\to H$ such that $\pr_+ T - 1\in L_2$ 
and $\pr_- T\in L_4$.  We equip this subspace with the 
topology coming from the metric  
$$ 
||T - T'|| = ||\pr_+ T - \pr_+ T'||_2 + ||\pr_- T - \pr_- T'||_4 
$$ 
where $||\cdot ||_2$ and $|| \cdot ||_4$ denote the norms 
of the Banach spaces $L_2$ 
and $L_4$ respectively.  
$\St_2$ has a natural 
structure of a Banach manifold, since it 
is an open subset of a Banach space. 
The projection $\St_2\to\Gr_2$ sends an admissible frame 
$w$ to the orthogonal projection $P_W$ onto its image $W$ 
($P_W$ is identified with an involution in the usual way).

Suppose that $X\colon H_+\to H$ is a linear map such that 
$\pr_+X\in L_2$ and $\pr_- X\in L_4$.  If $w\in \St_2$ 
then $w + tX$ is an injective map if the real number $t$ is 
small enough.  Therefore we see that the tangent 
space to $\St_2$ at $w$ can be identified with the 
Banach space of all $X\in \cB(H_+,H)$ with $\pr_+ X\in L_2$ 
and $\pr_- X\in L_4$ described above.

There is a natural connection $1$-form $\Theta$ on the principal $GL^2$ bundle 
$\St_2\to \Gr_2$ defined by 
$$ 
\Theta = w^{-1}P_Wdw
$$ 
Clearly this is $\Ad$-invariant and restricts to the Maurer-Cartan form 
on each fibre.  In order to be a connection $1$-form however 
$\Theta$ needs to take values in the Hilbert-Schmidt operators 
on $H_+$.  To see that this is the case observe that 
$$ 
w^{-1}P_Wdw = (\pr_+w)^{-1}\pr_+P_Wdw + (\pr_- w)^{-1}\pr_-P_W dw.   
$$ 
The following Lemma  
is a standard calculation.   
\begin{lemma} 
The curvature $2$-form $\Omega$ of the connection $\Theta$ 
is given by 
$$ 
\Omega = w^{-1}P_WdP_WdP_Ww. 
$$ 
\end{lemma} 


\section{The regularized determinant line bundle $\DET_2$} 
\label{sec3}
If $g$ is a bounded operator which differs from the identity 
by a trace class operator then we may form its determinant 
$\det(g)$.  In fact this determinant operator restricts to a homomorphism of groups  
$$ 
\det\colon GL^1\to \CC^* 
$$ 
where the group $GL^1$ consists of invertible operators $g$ such that 
$g-1$ is traceclass.  If the norm $|A|$ of $A$ is sufficiently small then this determinant is defined by the usual 
formula 
$$ 
\det(1+A) = \exp\left(\tr\left( \sum^\infty_{i=1}\frac{(-1)^{i-1}}{i}A^i\right) 
\right) 
$$ 
for trace-class operators $A$.  This expression 
obviously has no meaning if $A\in L_2$.  However, 
one can regularize this determinant by 
removing the divergent part of the trace 
and defining 
$$ 
\det_2(1+A) = \exp\left(\tr\left(\sum^\infty_{i=2} \frac{(-1)^{i-1}}{i}A^i\right) \right) 
$$
for $A\in L_2$ close to zero.  In general, to define $\det_2$ we proceed 
as follows (see \cite{Simon}).  For any bounded operator $A$ and any positive integer $n$, 
define 
$$ 
\cR_n(A) = (1+A)\exp\left(\sum^{n-1}_{j=1} (-1)^j\frac{A^j}{j}\right) - 1 
$$  
It is shown in Lemma 9.1 of \cite{Simon} that if $A\in L_n$, then $\cR_n(A)\in L_1$.  One then 
defines \cite{Simon}, for any positive integer $n$ and $A\in L_n$, 
$$ 
\det_n(1+A) = \det(1 + \cR_n(A)).  
$$  
When $n=2$ this process defines a map 
$$ 
\det_2\colon GL^2\to \CC^*, 
$$ 
however this 
is not a homomorphism; instead we have (see \cite{Simon}) 
\begin{equation} 
\label{eq: mult property of det_2} 
\det_2(1+A)(1+B) = \det_2(1+A)\det_2(1+B) \exp(-\tr(AB)) 
\end{equation}
for $A,B\in L_2$.   
Even though the regularized determinant $\det_2\colon GL^2\to \CC^*$ is not a homomorphism, 
we can still use it to define a 
determinant line bundle $\DET_2$ on $\Gr_2$ as follows.  
Following \cite{MR} we define an action of $GL^2$ on the quotient 
$\St_2\times \CC^*$ 
by  
$$ 
(w,z)g = (wg,z\omega(w_+,g)^{-1}) 
$$ 
where $w_+ = \pr_+ w$.  Here $\omega(w_+,g)$ is the 
function defined 
by
$$    
\omega(w_+,g) = \det_2(g)\exp(-\tr(w_+ -1)(g-1)) 
$$
Then, as in \cite{MR}, we let $\DET_2$ denote the quotient space 
$\DET_2 = (\St_2\times \CC^*)/GL^2$.  It can be shown that $\DET_2$ is a smooth 
principal $\CC^*$ bundle on $\Gr_2$.  In fact, since $\Det_2$ is a quotient of 
two complex manifolds, it is a holomorphic line bundle.  

\section{A connection $1$-form on $\DET_2$}
\label{sec4}
A good way to define a connection $1$-form on 
$\DET_2$ is to regard $\DET_2$ as an example of what Michael Murray 
has called a \emph{pre-line bundle} \cite{CJMSW}.  A pre-line bundle 
on a manifold $M$ consists of a surjective submersion 
$\pi\colon Y\to M$ together with a smooth map 
$f\colon Y^{[2]}\to \CC^*$ satisfying the `cocycle condition' 
$$
f(y_2,y_3)f(y_1,y_3)^{-1}f(y_1,y_2) = 1
$$ 
for points $y_1$, $y_2$, $y_3$ all lying in the same fibre of 
$Y$ over $M$.  Here $Y^{[2]} = \{(y_1,y_2)|\ \pi(y_1) = \pi(y_2)\}$.  It is a smooth 
submanifold of $Y^2$.  Given a pre-line bundle 
$(Y,f)$ on $M$ we can construct a principal $\CC^*$ bundle 
(and hence an associated line bundle) by forming the 
product $Y\times \CC^*$ 
and introducing the equivalence relation which identifies 
$$
(y,z)\sim (y',zf(y,y')) 
$$ 
for $(y,y')\in Y^{[2]}$. The quotient of $Y\times \CC^*$ 
by this equivalence relation defines a principal 
$\CC^*$ bundle $P$ over $M$, and in fact every 
principal $\CC^*$ bundle on $M$ arises in this way.  One can construct a connection 
$1$-form on this bundle as follows.  Suppose that 
there is a $1$-form 
$A$ on $Y$ such that 
\begin{equation} 
\label{eq: pre conn} 
f^{-1}df = \pi_2^*A - \pi_1^*A 
\end{equation}
where $\pi_1,\pi_2\colon Y^{[2]}\to Y$ denote the maps which 
omit the first and second factors in $Y^{[2]}$ respectively 
(such a $1$-form $A$ will exist  
if $M$ admits partitions of unity).  
Given a $1$-form $A$ satisfying~\eqref{eq: pre conn}, the $1$-form 
\begin{equation}
\label{eq: conn}
A + z^{-1}dz 
\end{equation} 
on $Y\times \CC^*$ descends to the quotient and defines a 
connection $1$-form on the principal $\CC^*$ bundle 
$P$.  Since pre-line bundles are not a familiar notion we will 
include the details of this.  We will show that the 1-form~\eqref{eq: conn} 
descends to a 1-form on $P$.  Suppose  
that  $(y_1,z_1)$ and $(y_2,z_2)$ are 
points lying in the same fiber of $Y\times \CC^*$ over $P$, and 
$(Y_1,\a_1)$, $(Y_2,\a_2)$ are tangent vectors at $(y_1,z_1)$, 
$(y_2,z_2)$ respectively which pushforward to the same tangent 
vector on $P$.  We need to show that 
$$ 
A(Y_1) + \a_1 = A(Y_2) + \a_2.  
$$ 
Since $(y_1,z_1)$ and $(y_2,z_2)$ lie in the same fiber over $P$ we must have 
$\pi(y_1) = \pi(y_2)$ and $z_2 = z_1f(y_1,y_2)$.  Similarly we must have 
$\a_2 = \a_1 + f^{-1}df(Y_1,Y_2)$.  Therefore 
$$ 
A(Y_2) + \a_2  = A(Y_2) + \a_1 + f^{-1}df(Y_1,Y_2)  
 = A(Y_1) + \a_1 
$$ 
as required.  It is easy to show that push forward of~\eqref{eq: conn} is invariant 
and restricts to the Maurer-Cartan 1-form on the fibers.  

As remarked above, the principal $\CC^*$ bundle $\DET_2$ 
is an example of a pre-line bundle for the submersion 
$\St_2\to \Gr_2$.  The cocycle $f$ in this case is defined to 
be 
$$ 
f(w_1,w_2) = \omega( (w_1)_+,g)^{-1} 
$$ 
for $(w_1,w_2)\in \St_2^{[2]}$ and 
where $g$ is the unique element of $GL^2$ such that 
$w_2 = w_1g$.  
The cocycle condition $f(w_2,w_3)f(w_1,w_3)^{-1}
f(w_1,w_2) = 1$ is easy to check.  
As mentioned in the introduction, it is possible to define a Hermitian structure on $\Det_2$ 
(see \cite{MR}).  As a Hermitian holomorphic line bundle $\Det_2$ therefore has a canonical 
connection, however obtaining a closed formula for its curvature seems to be rather difficult.  
Therefore we have constructed a connection on $\Det_2$ using the theory of pre-line bundles, as 
we now explain.  

We need to calculate $f^{-1}df$.  This is a somewhat longwinded calculation, 
and for this reason we have relegated it to Appendix B.  
The result is that 
$f^{-1}df$ is equal to 
\begin{equation} 
\label{eq: derivative of f} 
\tr((w_2)_+ - 1)\pi_1^*\Theta  - ((w_1)_+ - 1)\pi_2^*\Theta - \pr_+ \pi_2^*\pr_{H_w}dw + \pr_+\pi_1^*\pr_{H_w}dw) 
\end{equation}
Here $\Theta$ is an arbitrary connection 1-form on the principal bundle $\St_2\to \Gr_2$, and 
$\pr_{H_w}dw$ denotes the operator valued 1-form on $\St_2$ defined by the orthogonal projection of $dw\colon T\, \St_2\to \cB(H_+,H)$ 
onto the horizontal subspace (with respect to $\Theta$) at $w$. 
In order to write this as $\pi_2^*A - \pi_1^*A$ for some $1$-form $A$ on 
$\St_2$ we need a connection $\Theta$ on the principal $GL^2$ bundle 
$\St_2\to \Gr_2$ such that $\pr_+\pr_{H_w}dw$ is trace class.  Such a 
connection is given for example by 
\begin{equation} 
\label{eq: connection on St_2} 
\Theta = w^{-1}P_Wdw - w^{-1}P_W\pr_+dP_{W^\perp}w = w^{-1}P_Wdw + w^{-1}P_W\pr_+P_{W^\perp}dw. 
\end{equation}
To see this note that first of all $\pr_+dP_{W^\perp}w = \pr_+dP_{W^\perp}\pr_+w + 
\pr_+dP_{W^\perp}\pr_-w$ takes values in $L_2$, since $\pr_+dP_{W^\perp}\pr_+ \in L_2$ 
and $\pr_+dP_{W^\perp}\pr_-, \pr_- w\in L_4$.  Also, for this choice of 
$\Theta$ note that if $X$ is a tangent vector at $w$ then the horizontal 
projection $\pr_{H_w}X$ onto $H_w$ is given by 
$$ 
X - ww^{-1}P_WX - ww^{-1}P_W\pr_+P_{W^\perp}X = P_{W^\perp}X - P_W\pr_+P_{W^\perp}X. 
$$ 
Therefore 
\begin{align*}
\pr_+\pr_{H_w}X & = \pr_+P_{W^\perp}X - \pr_+P_W\pr_+P_{W^\perp}X \\ 
 & = \pr_+P_{W^\perp}\pr_+P_{W^\perp}X \\ 
 & = \left(\frac{1 - F_{11}}{2}\right) \left(\frac{1 - F_{11}}{2} 
 \right)\pr_+X + \left(\frac{1- F_{11}}{2}\right) 
 \left(\frac{ - F_{12}}{2}\right)\pr_- X 
 \end{align*} 
Since $F_{11} - 1\in L_2$ we see that this takes trace class values.  Finally then 
we can write down a connection $1$-form on $\Det_2$.   We describe this connection 1-form in the following 
proposition, where we also give an explicit formula for its curvature.\footnote{Jouko Mickelsson has informed me 
that he knew this expression for the first Chern class of $\Gr_2$, in a slightly different, but equivalent form.} 
\begin{proposition} 
\label{prop: connection}
A connection $1$-form on the principal $\CC^*$-bundle $\Det_2\to \Gr_2$ 
is given by 
$$ 
-\tr(\pr_+dw  - w^{-1}P_Wdw  - w^{-1}P_W\pr_+P_{W^\perp}dw) + z^{-1}dz.
$$ 
The curvature of this connection $1$-form is the $2$-form on 
$\Gr_2$ defined by 
$$ 
-\frac{1}{16}\tr((F - \epsilon)^2FdFdF) 
$$ 
\end{proposition} 

From equation~\eqref{eq: derivative of f} we see that for the choice~\eqref{eq: connection on St_2} of the connection $\Theta$ on 
$\St_2$, the following 1-form $A$ on $\St_2\times \CC^*$ satisfies the equation $\pi_2^*A - \pi_1^*A = f^{-1}df$: 
$$ 
A = -\tr((w_+ - 1)\Theta + \pr_+ \pr_{H_w} dw).  
$$ 
Therefore, 
by the general principles of pre-line bundles described above a connection $1$-form 
for $\Det_2$ is given by the push forward of the 1-form 
$$ 
-\tr((w_+ - 1)\Theta + \pr_+ \pr_{H_w}dw) + z^{-1}dz. 
$$ 
It is an easy calculation, using the definition of $\Theta$, to see that this is equal to 
\begin{multline*} 
-\tr\left((w_+ - 1)(w^{-1}P_Wdw + w^{-1}P_W\pr_+P_{W^\perp}dw) \right. \\ 
\left. + \pr_+P_{W^\perp}\pr_+P_{W^\perp}dw\right) + z^{-1}dz 
\end{multline*}
By straightforward manipulations one can show that this expression 
is the same as the one in the Proposition.  To find the curvature we 
need to find $d$ of the $1$-form $-\tr(\pr_+dw - w^{-1}P_Wdw - w^{-1}P_W\pr_+P_{W^{\perp}}dw)$.  
In order to do this we will make use of the fact that $ww^{-1}P_W = P_W$, 
which gives on differentiation $d(w^{-1}P_W) = w^{-1}P_WdP_W - w^{-1}P_Wdww^{-1}P_W$.  We 
calculate 
\begin{multline*} 
 -\tr(-w^{-1}P_WdP_Wdw + w^{-1}P_Wdww^{-1}P_Wdw 
 - w^{-1}P_WdP_W\pr_+P_{W^\perp}dw \\
+ w^{-1}P_Wdww^{-1}P_W\pr_+P_{W^\perp}dw + w^{-1}P_W\pr_+dP_Wdw) 
\end{multline*} 
since $dP_{W^\perp} = - dP_W$.  The term $w^{-1}P_Wdww^{-1}P_Wdw$ is trace 
class and its trace is easily seen to vanish due to the following property\footnote{  More generally (see \cite{Simon} Corollary 3.8) if $A,B\in \cB(H)$ 
have the property that $AB\in L_1$ and $BA\in L_1$, then $\tr(AB) = \tr(BA)$. } of 
the operator trace: if $A$ and $B$ are Hilbert-Schmidt operators 
then $\tr(AB) = \tr(BA)$ Also, we can differentiate the identity 
$P_W w = w$ to obtain $dw = dP_Ww + P_Wdw$, and using this we 
can re-write the term $w^{-1}P_WdP_Wdw$ as 
\begin{align*} 
w^{-1}P_WdP_Wdw & = w^{-1}P_WdP_WdP_W w + 
w^{-1}P_WdP_WP_Wdw \\ 
& = w^{-1}P_WdP_WdP_Ww,    
\end{align*} 
where we have used the fact that $P_WdP_WP_W = 0$.  
Finally, we can observe that the 
term $w^{-1}P_Wdww^{-1}P_W\pr_+ 
P_{W^\perp}dw = w^{-1}P_Wdww^{-1}P_W\pr_+dP_Ww$ 
is trace class and, using the cyclic property of the trace mentioned above,  
we can write   
$$ 
\tr(w^{-1}P_Wdww^{-1}P_W\pr_+dP_Ww) 
= \tr(w^{-1}P_W\pr_+dP_WP_Wdw). 
$$ 
Therefore our expression for the curvature becomes 
\begin{multline*} 
-\tr(-w^{-1}P_WdP_WdP_Ww - 
w^{-1}P_WdP_W\pr_+P_{W^\perp}dw \\ 
- w^{-1}P_W\pr_+dP_WP_Wdw + 
w^{-1}P_W\pr_+dP_Wdw) 
\end{multline*} 
Using the identity $dw = dP_Ww + P_Wdw$ we can simplify the last 
two terms in the above expression to $w^{-1}P_W\pr_+dP_WdP_Ww$.  Thus 
our new expression for the curvature is 
$$ 
-\tr(-w^{-1}P_W\pr_-dP_WdP_Ww - w^{-1}P_WdP_W\pr_+dP_Ww). 
$$ 
Each of the terms inside the trace belongs to the trace class ideal.  
Since $dP_WdP_WP_W = P_WdP_WdP_W$ and $dP_WP_WdP_W = P_{W^\perp}dP_WdP_W$ 
it is not hard to see that we can re-write the above expression 
as 
$$ 
-\tr(-\pr_-P_WdP_WdP_W +\pr_+P_{W^\perp}dP_WdP_W) 
$$
Again, each of the expressions $\pr_-P_WdP_WdP_W$ and 
$\pr_+P_{W^\perp}dP_WdP_W$ is trace class and so we 
may finally re-write this expression as 
$$ 
-\tr(-P_W\pr_+P_WdP_WdP_W + P_{W^\perp}\pr_-P_{W^\perp}dP_WdP_W) 
$$ 
In terms of the involutions $F$ associated to the 
projections $P_W$, this expression becomes 
$$ 
-\frac{1}{16}\tr((F-\epsilon)^2FdFdF). 
$$ 
The extra factors of $F-\epsilon$ serve to regularize 
the trace of the usual curvature $2$-form $FdFdF$ 
of the finite dimensional Grassmannian.  
As a consistency check one can also see that this gives 
a closed form. 

\section{The Mickelsson-Rajeev cocycle} 
\label{sec5}
In \cite{MR} the Lie algebra $2$-cocycle $\omega_A$ associated to the abelian extension 
$1\to \Map(\Gr_2,\CC^*)\to \widehat{GL}_{(2)}\to GL_{(2)}\to 1$ 
was computed. 
This is the well-known \emph{Mickelsson-Rajeev 
cocycle} 
$$ 
\omega_A = \frac{1}{8}\tr_C[[\epsilon,X],[\epsilon,Y]] (\epsilon - F)
$$ 
where $\tr_C$ denotes the \emph{conditional trace}.  Recall that 
$\tr_C$ is a regularization of the ordinary operator trace $\tr$ 
defined for operators $A$ on a polarized Hilbert space $H = H_+\oplus H_-$ by  
$$ 
\tr_C(A) = \frac{1}{2}\tr(A + \epsilon A\epsilon)
$$ 
whenever the latter trace exists.  
An operator $A$ for which $\tr_C(A)$ is defined is called conditionally 
trace class.  Clearly every trace class operator $A$ is conditionally 
trace class and moreover $\tr(A) = \tr_C(A)$ in this case. 
It is an easy calculation to show that $\omega_A$ may be re-written 
in terms of the usual operator trace as 
\begin{equation} 
\label{eq: cocycle 1}
\omega_A = \frac{1}{16}\tr(F- \epsilon )^2\epsilon[[\epsilon,X],[
\epsilon,Y]] 
\end{equation}
We will now compare $\omega_A$ to the cocycle $\omega_G$  
defined using the connection described in the preceding 
section.  In order to do this we need to compute the 
fundamental vector field $\hat{X}$ on $\Gr_2$ associated 
to the infinitesimal action of a vector $X\in \mathfrak{g}_{(2)}$.   
Since $g\in GL_{(2)}$ acts on $F\in \Gr_2$ by $g(F) = gFg^{-1}$ we 
see that $\hat{X}_F = [F,X]$.  Therefore the cocycle $\omega_G$ 
is given by 
\begin{equation} 
\label{eq: cocycle 2}
\omega_G = \frac{1}{16}\tr(F-\epsilon)^2F[[F,X],[F,Y]] 
\end{equation}
Although there are similarities between the two cocycles~\eqref{eq: cocycle 1} and~\eqref{eq: cocycle 2}, it seems as least as 
hard to guess a coboundary relating them as to use   
the formula~\eqref{eq: coboundary formula} of Proposition~\ref{prop: coboundary}.  
In Appendix A we use the latter method to derive  
the following expression for a coboundary $b$ relating the two cocycles: 
\begin{equation} 
b(X)(F) = \frac{1}{16}\tr((F-\epsilon)^3(F+\epsilon)\epsilon[\epsilon,X]) - \frac{1}{16}\tr((F-\epsilon)^4\epsilon X). 
\end{equation}

\appendix 

\section{Calculating the coboundary}

In this section we compute the coboundary 
$$ 
b(X) = A(\widehat{d\sigma(X)}) 
$$ 
of Proposition~\ref{prop: coboundary} for the Mickelsson-Rajeev extension~\eqref{eq: inf MR extn} when $\Det_2$ is 
equipped with the connection 1-form defined in Proposition~\ref{prop: connection}.  For this we 
need to review some constructions from \cite{MR}, in particular we need to understand the local structure of the 
group $\widehat{GL}_{(2)}$.  Firstly, 
we denote by $\cE_2$ the group 
$$ 
\cE_2 = \{(g,q)|\ g\in GL_{(2)},\ q\in GL(H_+), aq^{-1} - 1\in L_2\} 
$$ 
where $g$ is written in block diagonal form as 
$$ 
g = \left(\begin{array}{cc} a & b \\ c & d\end{array}\right).
$$  
The product $\cE_2\times \Map(\Gr_2,\CC^*)$ acts on the left of $\Det_2$ through 
the formula 
$$ 
(g,q,f)\cdot (w,\lambda) = (gwq^{-1}, f(F)\lambda \a(g,q;w)) 
$$ 
where $F$ is the Hermitian involution corresponding to the plane 
$W$ spanned by the admissible basis $w$, and where $\a(g,q;w)$ 
is a function satisfying a certain equivariance property.  The function $\a$ is not unique, and Mickelsson and Rajeev make the choice 
$$ 
\a = \exp\left(-\tr((1-q^{-1}a)(w_+ - 1) + q^{-1}b(\frac{1}{2}F_{21} - w_-))\right) 
$$ 
The group $\widehat{GL}_{(2)}$ is constructed from the product 
$\cE_2\times \Map(\Gr_2,\CC^*)$ as a quotient: 
$$ 
(\cE_2\times \Map(\Gr_2,\CC^*))/N, 
$$ 
where $N$ is a normal subgroup.  The precise description of $N$ will 
not be necessary for our calculation, and we refer to \cite{MR} 
for more details.   In an open neighbourhood of the identity element 
of $GL_{(2)}$, Mickelsson and Rajeev define a local section $\sigma\colon 
GL_{(2)}\to \widehat{GL}_{(2)}$ by 
$$ 
\sigma(g)= (g,a,1) \mod N.   
$$ 
Before we get to the calculation there is one last construction we need to review from \cite{MR}.  
Suppose that $W$ is a subspace of $H$ representing a point of $\Gr_2$ and that $w\colon H_+\to H$ is an 
admissible frame for $W$.  Then we can define (see \cite{MR, PS}) an invertible 
operator 
$$ 
h = \left(\begin{array}{cc} w_+ & \a \\ w_- & \b \end{array}\right) 
$$ 
such that $h(H_+) = W$ and $h(H_-) = W^\perp$.  If we denote 
$$ 
h^{-1} = \left(\begin{array}{cc} x & y \\ u & v \end{array}\right) 
$$ 
so that $x = w^{-1}P_W\pr_+$ and $y = w^{-1}P_W\pr_-$, 
then the involution $F$ corresponding to $W$ can be written $F = h\epsilon h^{-1}$.  
As noted in \cite{MR} this implies the pair of equations $F_{11} = 2w_+x - 1$ and 
$F_{21} = 2w_-x$.  

We need to compute the fundamental vector field $\widehat{d\sigma(X)}$ 
where $d\sigma$ denotes the derivative of $\sigma$ and $X$ is a vector 
in $\mathfrak{gl}_{(2)}$.  An easy computation gives that 
$$ 
\widehat{d\sigma(X)} = (Xw - wX_{11}, -\tr(X_{12}(\frac{1}{2}F_{21} - w_-))). 
$$ 
We now need to compute the result of applying the connection 1-form $A$ 
on $\Det_2$ to this vector field.  Note that $A$ can be written as  
$$ 
-\tr((1-  w^{-1}\pr_W\pr_+)\pr_+dw) +  \tr(w^{-1}P_W\pr_-P_Wdw) + z^{-1}dz   
$$ 
We first compute $\tr(w^{-1}P_W\pr_-P_W(Xw - wX_{11}))$.  We can write this as 
\begin{equation} 
\label{eq: expression 1}
\tr(w^{-1}P_W\pr_-P_WX\pr_+w_+ + w^{-1}P_W\pr_-P_WX\pr_-w_-  - w^{-1}P_W\pr_-wX_{11}) 
\end{equation}
and then split~\eqref{eq: expression 1} up into the following sum of traces (it is straightforward to 
check that all of the expressions involved are of trace class): 
\begin{multline} 
\label{eq: expression 2}
\tr(w^{-1}P_W\pr_-P_W\pr_-X_{21}w_+) + \tr(w^{-1}P_W\pr_-P_WX\pr_-w_-) \\
+\tr(w^{-1}P_W\pr_-P_W\pr_+X_{11} - w^{-1}P_W\pr_-wX_{11})  \\
+ \tr(w^{-1}P_W\pr_-P_W\pr_+X_{11}(w_+ - 1)) 
\end{multline}
By using the cyclic property of the trace we can further write~\eqref{eq: expression 2} as 
\begin{multline}
\label{eq: expression 3}
\tr(\pr_+P_W\pr_-P_W\pr_-X_{21} + \pr_-P_W\pr_-P_WX\pr_- \\ 
+ (w_+ - 1)w^{-1}P_W\pr_-P_W\pr_+X_{11}) \\
+ \tr(w^{-1}P_W\pr_-P_W\pr_+X_{11} - w^{-1}P_W\pr_-wX_{11}) 
\end{multline} 
The term $\tr((w_+ - 1)w^{-1}P_W\pr_-P_W\pr_+X_{11})$ in~\eqref{eq: expression 3} can be written as 
$$ 
\tr(\pr_+P_W\pr_-P_W\pr_+X_{11} - w^{-1}P_W\pr_-P_W\pr_+X_{11}). 
$$ 
This can be combined with the last term in~\eqref{eq: expression 3}
to get 
$$ 
\tr(\pr_+P_W\pr_-P_W\pr_+X_{11} - w^{-1}P_W\pr_-wX_{11}) = \tr( \frac{F_{12}F_{21}}{4}- yw_- )X_{11}). 
$$ 
Therefore we find that we can write~\eqref{eq: expression 3} as  
\begin{equation} 
\label{eq: expression 4}
\tr(\pr_+P_W\pr_-P_W\pr_-X_{21} + \pr_-P_W\pr_-P_WX\pr_- + (\frac{F_{12}F_{21}}{4}- yw_- )X_{11}).
\end{equation}
We now compute the contribution to $A(\widehat{d\sigma(X)})$ coming from 
\begin{equation} 
\label{eq: expression 5}
-\tr((1-w^{-1}P_W\pr_+)(\pr_+Xw - w_+X_{11})) 
\end{equation}
We write $x$ for $w^{-1}P_W\pr_+$.  From $w_+x = \pr_+P_W\pr_+$ we see that $x \in 1 + L_2$.  
We also write $\pr_+Xw = X_{11}w_+ + X_{12}w_-$.  Then, after a little 
manipulation,~\eqref{eq: expression 5} can be written as 
\begin{equation} 
\label{eq: expression 6} 
- \tr\left((xw_+ - w_+x)X_{11}\right) - \tr\left((1-x)X_{12}w_-\right)
\end{equation} 
Since $F_{21} = 2w_- x$ we see 
that we can write the term $\tr(X_{12}(F_{21}/2 - w_-))$ as 
$$ 
\tr(X_{12}(F_{21}/2 - w_-)) = - \tr(X_{12}w_-(1-x)) = - \tr((1-x)X_{12}w_-). 
$$ 
Therefore, on applying $A$ to the fundamental vector field $\widehat{d\sigma(X)}$ we 
see that the contribution coming from $z^{-1}dz$ cancels with the term $-\tr((1-x)X_{12}w_-)$ 
appearing in the expression for $-\tr(1-w^{-1}P_W\pr_+)(\pr_+Xw - w_+X_{11})$.  We 
obtain the following expression for $A(\widehat{d\sigma(X)})$: 
\begin{multline*}
\tr(\pr_+P_W\pr_-P_W\pr_-X_{21} + \pr_-P_W\pr_-P_WX\pr_-) \\
+ \tr((w_+x - xw_+ -yw_- + \frac{F_{12}F_{21}}{4})X_{11}) 
\end{multline*} 
Since $xw_+ + yw_- = 1$, the last term can be written as 
$$ 
\tr((w_+x - 1 + \frac{F_{12}F_{21}}{4})X_{11}) 
=- \frac{1}{4}\tr((1-F_{11})^2X_{11}), 
$$ 
where we have used the identity $w_+x = (1+F_{11})/2$.  
Therefore we have the following expression for the coboundary $b(X) = A(\widehat{d\sigma(X)})$, independent 
of $w$: 
$$ 
\tr(\pr_+P_W\pr_-P_W\pr_-X_{21} + \pr_-P_W\pr_-P_WX\pr_-) 
- \frac{1}{4}\tr((1-F_{11})^2X_{11}). 
$$
We can further write the term $\tr(\pr_-P_W\pr_-P_WX\pr_-)$ as 
$$ 
\tr(\pr_-P_W\pr_-P_W\pr_+X_{12}) + \tr(\pr_-P_W\pr_-P_W\pr_-X). 
$$ 
Since $\pr_-P_W\pr_- = (F_{22} + 1)/2$ we can rewrite the formula for 
$b(X)$ as 
\begin{multline*} 
\tr(\pr_+P_W\pr_-P_W\pr_-X_{21} + \pr_-P_W\pr_-P_W\pr_+X_{12}) \\ 
+ \frac{1}{4}\tr((F_{22} + 1)^2X) - \frac{1}{4}\tr((F_{11} - 1)^2X). 
\end{multline*}
We can write $b(X)$ in terms of $F$, $\epsilon$ and $X$ by setting $\pr_+ = \frac{1+\epsilon}{2}$, $P_W = \frac{1+F}{2}$ etc.  
For instance 
$$
F_{11} - 1= \pr_+(F-1)\pr_+ = \frac{1}{4}(1+\epsilon)(F-1)(1+\epsilon) = -\frac{1}{4}(F-\epsilon)^2(1+\epsilon) 
$$ 
using the identities $(F-\epsilon)^2 = 2 - \epsilon F - F\epsilon$ and $\epsilon F\epsilon = 2\epsilon - F - (F-\epsilon)^2\epsilon$.  
Similarly $F_{22} + 1 = (F-\epsilon)^2(1-\epsilon)/4$.  Thus the term $ \frac{1}{4}\tr((F_{22} + 1)^2X) - \frac{1}{4}\tr((F_{11} - 1)^2X)$ 
above can be written as 
$$ 
-\frac{1}{16}\tr((F-\epsilon)^4\epsilon X).
$$ 
Since $\pr_-X_{21} = \pr_-[\epsilon,X]$ and $\pr_+P_W\pr_- = \frac{1}{8}(1+\epsilon)F(1-\epsilon)$ we can 
write 
\begin{align*} 
& \tr(\pr_+P_W\pr_-P_W\pr_-X_{21}) \\ 
= & \frac{1}{32}\tr((1+\epsilon)F(F-\epsilon)^2(1-\epsilon)[\epsilon,X]) \\ 
= & \frac{1}{32}\tr((F-\epsilon)^2(1+\epsilon)F(1-\epsilon)[\epsilon,X]) \\ 
= & \frac{1}{32}\tr((F-\epsilon)^3(F+\epsilon)(\epsilon - 1)[\epsilon,X]). 
\end{align*} 
We obtain a similar expression for $\tr(\pr_-P_W\pr_-P_W\pr_+X_{12})$: 
$$ 
\tr(\pr_-P_W\pr_-P_W\pr_+X_{12}) = \frac{1}{32}\tr((F-\epsilon)^3(F+\epsilon)(\epsilon +1)[\epsilon,X]). 
$$ 
Combining all of these expressions we can finally write 
$$ 
b(X) = \frac{1}{16}\tr((F-\epsilon)^3(F+\epsilon)\epsilon[\epsilon,X]) - \frac{1}{16}\tr((F-\epsilon)^4\epsilon X). 
$$ 

\section{Derivation of equation~\eqref{eq: derivative of f} } 

Recall that $f\colon \St_2^{[2]}\to GL^2$ is defined to be 
$$ 
f(w_1,w_2) = \omega((w_1)_+,g)^{-1} = (\det_2(g))^{-1}\exp(\tr((w_1)_+ - 1)(g-1)), 
$$ 
where $g = g(w_1,w_2)$ is the unique element of $GL^2$ such that $w_2 = w_1g(w_1,w_2)$.  
Therefore $f^{-1}df$ is equal to 
\begin{equation} 
\label{eq: expression B1}
-(\det_2(g))^{-1}d\det_2(g) + \tr \, \pi_2^*dw_+(g-1) + \tr((w_1)_+ - 1)dg 
\end{equation}
Thus we need to calculate the derivatives $d\det_2(g)$ and $dg$.  The latter is easy, it is 
well known that 
\begin{equation} 
\label{eq: expression B2} 
dg(X_1,X_2) = g\Theta(X_2) - \Theta(X_1)g
\end{equation}
where $(X_1,X_2)$ is a tangent vector at $(w_1,w_2)$ in $\St_2^{[2]}$ and $\Theta$ is an arbitrary 
connection 1-form on the principal $GL^2$ bundle $\St_2\to \Gr_2$.  
To calculate $d\det_2(g)$ we need to calculate the derivative 
$$ 
\frac{d}{dt}\Big|_{t=0} \det_2(\exp(-t\Theta(X_1))g\exp(t\Theta(X_2)))). 
$$ 
For this we will use the multiplicative property~\eqref{eq: mult property of det_2} of the regularized 
determinant $\det_2$ to write $\det_2(\exp(-t\Theta(X_1))g\exp(t\Theta(X_2))))$ as the product  
\begin{multline*} 
\det_2(\exp(-t\Theta(X_1)))\cdot \det_2(g)\cdot \det_2(\exp(t\Theta(X_2)))\cdot \\ 
\cdot \exp\left( -\tr(g-1)(\exp(t\Theta(X_2)) - 1)\right)\cdot  \\
\cdot \exp\left( -\tr((\exp(-t\Theta(X_1)) - 1)(g\exp(t\Theta(X_2))-1)\right) 
\end{multline*} 
We recall that $\det_2(\exp(tA))$ is $O(t^2)$ and so on taking derivatives we end up with 
\begin{equation} 
\label{eq: expression B3} 
\det_2(g)\left( -\tr((g-1)\Theta(X_2)) + \tr(\Theta(X_1)(g-1))\right) 
\end{equation}
Combining~\eqref{eq: expression B1},~\eqref{eq: expression B2} and~\eqref{eq: expression B3} we 
get the following formula for $f^{-1}df$: 
\begin{multline} 
\label{eq: expression B4} 
\tr((g-1)\Theta(X_2)) - \tr(\Theta(X_1)(g-1)) + \tr\, \pi_2^*dw_+(g-1) \\ 
+ \tr((w_1)_+ - 1)(g\Theta(X_2) - \Theta(X_1)g) 
\end{multline} 
In this formula we can write $\pi_2^*dw_+ = (w_1)_+A(X_1) + \pi_2^*\pr_{H_w}dw$, where 
$\pr_{H_w}dw$ denotes the operator valued form on $\St_2$ which is the composition of the inclusion 
$dw\colon T\, \St_2\to \cB(H_+,H)$ followed by orthogonal projection onto the horizontal subspace 
$H_w$ defined by $\Theta$.  Then we see that~\eqref{eq: expression B4} can be re-written as 
\begin{multline}
\label{eq: expression B5} 
\tr(((w_2)_+ - 1)\pi_1^*\Theta  - ((w_1)_+ - 1)\pi_2^*\Theta - \pr_+ \pi_2^*\pr_{H_w}dw  \\ 
+\pi_2^*(dw_+) g - (w_1)_+ \pi_2^*\Theta g) 
\end{multline} 
Since horizontal subspaces are translation invariant, we can write $\pi_2^*(dw) g - \pi_2^*(w\Theta) g = \pi_1^*\pr_{H_w}dw$.  
Therefore~\eqref{eq: expression B5} becomes 
$$  
\tr((w_2)_+ - 1)\pi_1^*\Theta  - ((w_1)_+ - 1)\pi_2^*\Theta - \pr_+ \pi_2^*\pr_{H_w}dw + \pr_+\pi_1^*\pr_{H_w}dw) 
$$
which is equation~\eqref{eq: derivative of f}.  
  
\section*{Acknowledgments} 

I thank Michael Murray for some useful conversations and Alan Carey for his comments on 
a first draft of this note.  I also thank Christoph Schweigert for his excellent help in improving the 
exposition.  I am especially grateful to Jouko Mickelsson for some useful conversations and 
his detailed comments and suggestions on reading a previous version of this note. 
This work was supported by the Collaborative Research Center 
676 `Particles, Strings and the Early Universe'.

\end{document}